\documentclass[11pt,a4paper]{article}
\usepackage[margin=28mm]{geometry}
\usepackage[authoryear,round]{natbib}
\usepackage[T1]{fontenc}
\usepackage{amsmath,amssymb,amsthm,mathtools,booktabs,array}
\usepackage[hidelinks]{hyperref}

\newtheorem{theorem}{Theorem}[section]
\newtheorem{lemma}[theorem]{Lemma}
\newtheorem{proposition}[theorem]{Proposition}
\newtheorem{corollary}[theorem]{Corollary}
\theoremstyle{remark}\newtheorem{remark}[theorem]{Remark}
\newcommand{\E}{\mathbb E}
\newcommand{\Pp}{\mathbb P}

\newcommand{\Beta}{\operatorname{Beta}}

\newcommand{\norm}[1]{\lVert#1\rVert}

\title{The minimax value of outcome feedback in budget-constrained adaptive sampling}
\author{Shengye Tao\\[2pt]\small Beijing Institute of Mathematical Sciences and Applications, Beijing, China\\\small\texttt{shengyetao@bimsa.cn}}
\date{}
\begin{document}
\maketitle
\begin{abstract}
Whether adaptivity improves estimation is a classical question in sequential design and information-based complexity. In practice, adaptive data collection can use two different kinds of information: process data on whether each contact succeeded, and the observed outcomes themselves. How much each kind is worth in the worst case under a hard budget has not been settled. We study this question for estimating a common Bernoulli mean through channels with unknown response probabilities and known costs, where failed contacts are charged and the budget must hold on every sample path. We compare static designs, designs that adapt to response indicators only, and fully adaptive designs, each with an arbitrary estimator. Our main theorem shows that when the costs are invariant under a group of channel permutations acting transitively on the nuisance parameter, outcome-adaptive allocation has zero minimax value: response-adaptive and fully adaptive designs have the same minimax risk, given by a one-dimensional optimization and attained by an explicit least-favourable prior and equalizer rule. For two contacts we obtain a closed-form criterion for when response feedback strictly improves on every static design. An unequal-cost example with a computer-certified exact saddle shows that the same conclusion can hold without symmetry, even though no product prior is least favourable there. A correlated prior shows that outcome feedback can still reduce Bayes risk, so the equivalence is genuinely minimax. In this idealized setting, the results support responsive survey designs that adapt to paradata but not to outcomes.
\end{abstract}


\section{Introduction}
Does adaptivity help? The question is as old as sequential analysis, and it has received sharply different answers. In information-based complexity, adaptively chosen information is known to be no more powerful than nonadaptive information for linear problems in the worst case, and the same holds under Gaussian noise, whereas adaptation can be much more powerful when those conditions fail \citep{kadane1988}. In sequential experimental design and bandit problems, by contrast, adaptive allocation is the rule rather than the exception. A precise answer therefore depends on the information structure, on the loss and on the criterion.

This paper studies the question in a setting where two different kinds of feedback are available. Statistical observations are often obtained only after a successful contact, and different access channels---survey modes, data vendors or measurement devices---may measure the same outcome distribution while differing in availability and cost. When the budget is finite, a failed attempt consumes resources and changes what can still be observed. An adaptive design can use the \emph{response history} to learn which channel is productive, and it can also use the \emph{observed outcomes} when deciding where, or whether, to sample next. Responsive survey designs use process data of the first kind to guide mid-survey decisions \citep{groves2006}. Whether outcome-dependent allocation adds anything beyond such process data is a natural question, and for estimation of a common mean it has, to our knowledge, not been answered in exact form.

We consider squared-error estimation of a common Bernoulli mean under a hard, pathwise budget in which failures are charged at the same contact cost as responses. Three nested classes of designs are compared: designs fixed before any observation, allocations that use only channel labels and response indicators, and allocations that additionally use outcome values. Every class permits arbitrary randomization and arbitrary terminal estimation, and the response probabilities are fixed unknown parameters. Our contributions are as follows.
\begin{enumerate}
\item \emph{Bayes reduction} (Proposition~\ref{prop:bayes}). Under an independent Beta prior on the mean and an arbitrary prior on the response probabilities, outcome-dependent allocation is never needed for Bayes optimality, with arbitrary costs.
\item \emph{An exact minimax saddle} (Theorem~\ref{thm:orbit}). When a finite group of cost-preserving channel permutations acts transitively on the nuisance parameter, response-adaptive and fully adaptive designs have the same minimax risk, given by a one-dimensional maximization. We construct a least-favourable prior and an equalizer that mixes at most two response trees.
\item \emph{A sharp criterion for response feedback} (Corollary~\ref{cor:two}). With two equal-cost contacts and swapped response rates $(a,b)$, response feedback strictly improves on every static design if and only if $a\ne b$ and $4(a+b)-5ab>2$.
\item \emph{Beyond symmetry} (Theorem~\ref{thm:asym} and Proposition~\ref{prop:nonproduct}). In an unequal-cost experiment with budget five, a computer-certified exact saddle shows that outcome feedback again has zero minimax value, although every product prior falls strictly short of the least-favourable value.
\item \emph{Minimax versus Bayes} (Proposition~\ref{prop:counter}). A correlated prior under which outcome feedback strictly lowers Bayes risk shows that the equivalence is specific to the minimax criterion and is not a comparison of experiments.
\end{enumerate}
The proofs combine a posterior-variance factorization with group averaging and a slope-balancing argument that equalizes the center and endpoint risks of at most two active response trees. The resulting lower bound applies to every fully adaptive policy and every estimator.

The paper is organized as follows. Section~\ref{sec:model} sets up the model and the information classes. Section~\ref{sec:bayes} proves the Bayes reduction, Section~\ref{sec:orbit} the minimax saddle under symmetry, and Section~\ref{sec:two} the two-contact classification. Section~\ref{sec:asym} treats the unequal-cost experiment and Section~\ref{sec:counter} the Bayes counterexample. Section~\ref{sec:discussion} discusses scope and limitations.

\subsection{Relation to existing work}
Sequential allocation and the optimization of posterior terminal loss are established methods. \citet{hardwick1998common} study sequential choice between Bernoulli$(p)$ and Bernoulli$(p^r)$ experiments, including asymptotic Bayes results and finite comparisons using maximum likelihood estimation. A related Bayesian formulation of these same two experiments is considered by \citet{jokielmagiera2002}. \citet{hardwick1998missing} describe Bayes allocation with censored outcomes and a general dynamic-programming approach. Their examples have distinct outcome parameters and a shared censoring mechanism. Neither observation law is the common-mean, unknown-channel response experiment considered here, and a Bayes recursion alone does not give a matching frequentist saddle.

Cost-sensitive adaptive estimation of a common mean has also been studied from a learning perspective. \citet{neufeld2014} study bandit allocation among unbiased Monte Carlo estimators, with MSE-regret relative to the best estimator and an extension to stochastic observation costs. More recently, \citet{lin2026} analyze adaptive allocation for arm-level mean estimation and compare specified procedures with uniform sampling. Our question is complementary: rather than comparing specified procedures, we compute the exact minimax value of each information class over all estimators in one fixed hard-budget experiment.

Classical background includes ancillarity in multinomial estimation \citep{he1990} and minimax prediction for multinomial and multivariate hypergeometric laws \citep{jokiel1998}. Our random-count shrinkage calculation is proved directly after fixing a response-only tree. An outcome-dependent policy may make its terminal count depend on the mean, so it cannot automatically be reduced to that experiment. Finite-group risk averaging is also classical. \citet{bai2023} give finite-sample minimax results for permutation-invariant assignment models and distinguish randomization of a fixed estimator from joint design and linear estimation. We use risk averaging for controlled trees only when permutations preserve costs; the separate posterior reduction and equalizer construction are still needed.

Additional nearby work has a different inferential target. \citet{biswas2004} study adaptive clinical allocation with missing responses and imputation under missing at random; \citet{oluyede2003} compare sequential and fixed estimation of a difference between two population means, using Bayes risk and expected cost. These do not identify the hard-budget common-mean saddle or its response-versus-outcome feedback comparison. Adaptive allocation, missing-data recursions and conjugate variance identities are standard tools; our contribution lies in combining them into exact minimax statements.

The equalities below concern one minimized worst squared risk, not a stochastic simulation of experiments. The comparison in \citet{blackwell1953} is stronger: it concerns attainable risks for decision problems and stochastic transformations of experiments. Section~\ref{sec:counter} exhibits a joint prior with strictly smaller full-feedback Bayes risk, even in an experiment covered by our minimax theorem. Our equivalence is therefore specific to the minimax criterion.

\section{Statistical model and information classes}\label{sec:model}
Let $H$ be a positive integer, let $c_h>0$ and $0\le B<\infty$ be known, and let $\Gamma$ be a known nonempty Borel set. The parameter space is
\[
 \Theta=[0,1]\times\Gamma,\qquad \rho=(\rho_1,\ldots,\rho_H)\in\Gamma\subseteq[0,1]^H.
\]
Nature fixes $(p,\rho)$ before independent policy randomization. A fresh contact to channel $h$ produces $R\sim\mathrm{Bernoulli}(\rho_h)$ and $Y\sim\mathrm{Bernoulli}(p)$, independent conditional on the parameter and channel, and independent of past contacts given their selected channels. We observe $R$ and, when $R=1$, $Y$. Each contact costs $c_h$, including $R=0$. A policy chooses an affordable channel or stops, and every realized path must have total cost at most $B$. There are at most $T=\lfloor B/\min_h c_h\rfloor$ contacts.

The final estimate $d$ targets $p$ with loss $(d-p)^2$. It may use the complete legal observation history and the realized independent random seed. Projection onto $[0,1]$ does not increase loss, so we restrict estimates to that interval. Define $\Pi_{\mathrm{static}}$ to choose the entire finite schedule before data, including any preassigned stopping. Define $\Pi_{\mathrm{response}}$ to allocate and stop using channel labels, response indicators, remaining budget and randomization. Define $\Pi_{\mathrm{full}}$ to allow observed $Y$ values in those decisions as well. All three classes allow the same terminal information. For $\mathcal C\in\{\mathrm{static},\mathrm{response},\mathrm{full}\}$ put
\begin{equation}\label{eq:minimax}
 V_{\mathcal C}=\inf_{(\pi,d)\in\mathcal C}\sup_{(p,\rho)\in\Theta}
       \E_{p,\rho}^{\pi}(d-p)^2.
\end{equation}
Randomization is averaged before the supremum. In particular,
\[
 V_{\mathrm{full}}\le V_{\mathrm{response}}\le V_{\mathrm{static}}.
\]
We abbreviate these values by $V_F,V_R,V_S$. The components of $(p,\rho)$ are fixed unknown constants in the frequentist experiment; independence of them below is a property of a specified prior, not a restriction on $\Theta$.

\begin{lemma}[Randomized rules and finite histories]\label{lem:law}
Every measurable rule in these classes has a well-defined finite-horizon
probability law. Its risk is the integral, over its independent seed, of
the risks of deterministic feasible trees with clipped terminal estimates.
For any finite prior, posterior terminal optimization and backward
minimization over affordable contacts give a lower bound for the Bayes
risk of every such randomized rule. If the parameter is a Bernoulli mean
and a nuisance response vector, the risk of a fixed-budget full rule is,
for fixed nuisance, a polynomial in $p$ of degree at most $T+2$.
\end{lemma}
\begin{proof}
Let $U$ have the policy's parameter-independent seed law $\mu$.
Conditional on a fixed $u$, a rule selects an action at each finite
history. Each contact has the three observable outcomes $0,10,11$, with
probabilities $1-\rho_h$, $\rho_h(1-p)$ and $\rho_hp$.
The probability of a leaf is the product of its contact probabilities.
Induction on depth shows that leaf probabilities sum to one: a stop
preserves its mass and a contact splits its mass into three children
whose probabilities sum to one. All paths have depth at most $T$,
because each contact costs at least $\min_h c_h$.
For fixed $u$, this construction is a finite probability space;
measurability of decisions and terminal estimates makes its probabilities
and losses measurable functions of $u$. Integrating the finite kernel
against $\mu$ therefore gives the required law. The loss after clipping
is at most one, so all risk integrals exist and finite sums commute with
the seed integral.

For a finite prior, conditional likelihoods are finite sums and
nonzero-mass posteriors are ordinary ratios. At a leaf, square completion
shows that the posterior mean minimizes loss. At an earlier history,
every chosen contact has expected continuation loss at least the minimum
over the affordable contacts and stopping. Inducting backward, then
integrating over $u$, proves the assertion for the entire seed space;
no finite-mixture restriction is imposed. Zero-mass histories contribute
zero and can receive any legal continuation.

There are finitely many possible channel/observation histories of depth
at most $T$. For each, its seed-dependent inclusion indicator and
terminal estimate are parameter independent. Its likelihood is a fixed
nuisance factor times $p^s(1-p)^{n-s}$, with $n\le T$.
Multiplication by the squared estimation loss gives a polynomial of
degree at most $T+2$. Its coefficients are bounded measurable functions
of $u$; integrating their finite sum proves the last assertion.
\end{proof}

\section{Bayes feedback reduction}\label{sec:bayes}
Let $z$ denote the channel/response history, $n$ the number of responses and $s$ the number of observed successes. For a given full realized history, including its seed, the likelihood has the form
\begin{equation}\label{eq:factor}
 K(\text{history})p^s(1-p)^{n-s}
       \prod_h\rho_h^{u_h}(1-\rho_h)^{v_h},
\end{equation}
where $u_h,v_h$ count responses and failures and $K$ is parameter-free. Outcome-based choices can change the distribution of $n$, but not this conditional factorization.

\begin{proposition}[Product-prior Bayes reduction]\label{prop:bayes}
Suppose $p\sim\Beta(\alpha,\beta)$, $\alpha,\beta>0$, independently of an arbitrary nuisance prior $\nu$. Under the experiment in Section~\ref{sec:model}, the minimum squared-error Bayes risk over full-feedback procedures equals that over response-only procedures. Heterogeneous costs and multiple nuisance orbits are permitted. An attaining terminal estimator is
\[
 d=(\alpha+s)/(\alpha+\beta+n).
\]
\end{proposition}
\begin{proof}
By \eqref{eq:factor}, the posterior is the product of $\Beta(\alpha+s,\beta+n-s)$ and the response-history posterior $\nu_z$. Write $A=\alpha+s$, $D=\beta+n-s$, $t=A+D$ and $v=AD/[t^2(t+1)]$. A responding observation gives expected posterior variance
\begin{align}\label{eq:variance}
 \frac At\frac{(A+1)D}{(t+1)^2(t+2)}
 +\frac Dt\frac{A(D+1)}{(t+1)^2(t+2)}
  =v\frac{t}{t+1}.
\end{align}
A failure leaves the $p$ posterior unchanged. Put $\ell_h(z)=\E_{\nu_z}\rho_h$. A backward induction on the maximal number of remaining contacts shows that the conditional optimized loss is $vG(t,z,e)$, where $e$ is remaining budget and
\begin{equation}\label{eq:bayesbellman}
\begin{split}
 G(t,z,e)=\min\Bigl\{1,\min_{c_h\le e}
 \Bigl[&(1-\ell_h)G(t,zh0,e-c_h)\\
       &+\ell_h\frac{t}{t+1}G(t+1,zh1,e-c_h)\Bigr]\Bigr\}.
\end{split}
\end{equation}
Here $zh0,zh1$ append the response indicator. At a terminal state $G=1$. At a continuation, the induction multiplier after a response does not depend on its success value; \eqref{eq:variance} therefore averages it exactly. Costs and admissible actions depend on $e,z$, not $s$. Since $t=\alpha+\beta+n$ is response-history measurable, an action minimizing \eqref{eq:bayesbellman} is response-only. Randomized actions cannot improve its conditional minimum. The posterior mean minimizes terminal squared loss among all estimators. Impossible histories can be assigned any legal action. Positive costs and finite budget make the induction finite.
\end{proof}

The proposition asserts equality of optimized Bayes values for the specified product prior. It asserts neither minimax equality nor risk domination of each full policy by a response-only policy. A joint prior can make the nuisance posterior depend on $s$, invalidating the separation. Conversely, neither equal costs nor symmetry was used in the proof.

\section{A minimax saddle under cost-preserving symmetry}\label{sec:orbit}
Let $\mathsf G$ be a finite group of channel permutations with $c_{gh}=c_h$ for every $g,h$. It acts by $(g\rho)_h=\rho_{g^{-1}h}$. Assume that
\begin{equation}\label{eq:orbit}
 \Gamma=\{g\rho^0:g\in\mathsf G\}
\end{equation}
is a single orbit, and let $\nu$ be uniform on its distinct elements. Relabeling both the channels chosen by a policy and the channel labels in its history preserves admissibility and budget. If $D$ denotes a policy-estimator pair, its conjugate satisfies
\[
 r(D^g;p,\rho)=r(D;p,g^{-1}\rho).
\]
Draw $g$ uniformly before data and keep it in the known seed. This symmetrization averages risk, does not increase worst risk, and makes risk the same at every nuisance state in the orbit. It applies to full policies as well as response-only ones. It does not by itself make allocations outcome-independent.

Let $\mathcal T_R$ be the finite set of deterministic feasible response-only trees, including stopping, and let $\mathcal T_S$ be the finite set of static schedules. For each tree let $N$ denote its terminal number of responses. For $\alpha>0$ define
\begin{equation}\label{eq:bcount}
 b_\alpha(n)=\frac{\alpha}{2(2\alpha+1)(n+2\alpha)},\qquad
 F_C(\alpha)=\min_{\pi\in\mathcal T_C}\E_{\nu}^{\pi}b_\alpha(N),\quad C=R,S.
\end{equation}
At zero set $b_0(0)=1/4$, $b_0(n)=0$ for $n>0$. The expectation uses the nuisance prior and the response experiment only.

\begin{theorem}[Single-orbit minimax saddle]\label{thm:orbit}
Under \eqref{eq:orbit} and the cost-preserving group assumptions,
\begin{equation}\label{eq:mainorbit}
 V_F=V_R=\max_{\alpha\ge0}F_R(\alpha),\qquad
 V_S=\max_{\alpha\ge0}F_S(\alpha).
\end{equation}
The maxima are attained. At an interior maximizing $\alpha_*>0$, an attaining rule mixes at most two symmetrized active trees and uses
\[
 q_{\alpha_*}(n,s)=\frac{s+\alpha_*}{n+2\alpha_*},\qquad q_{\alpha_*}(0,0)=1/2.
\]
Its risk is constant in $p$ and nuisance. At a maximizing $\alpha_*=0$, an active symmetrized tree with endpoint risk no smaller than center risk attains the value, using $s/n$ at $n>0$ and $1/2$ at zero.
\end{theorem}

\paragraph{Idea of the proof}
The lower bound comes from priors of the form $\Beta(\alpha,\alpha)\times\nu$. By Proposition~\ref{prop:bayes}, outcome feedback is useless against such priors, so the best full-feedback Bayes risk is $F_R(\alpha)$ and $V_F\ge\max_\alpha F_R(\alpha)$. For the matching upper bound we need a response-only rule whose risk is constant both in the nuisance state and in $p$. Group averaging handles the nuisance state. In $p$, the risk of a shrinkage estimator after a response tree is a quadratic whose curvature is proportional to $E-A$, the gap between its endpoint and center risks, and by \eqref{eq:slope} this gap is proportional to the slope of the tree's Bayes risk in $\alpha$. At the maximizer of $F_R$, a mixture of at most two active trees has zero slope, hence a flat risk. The least-favourable prior thus sits exactly where the best response trees balance their center and endpoint risks.

\begin{proof}
For a response-only tree and fixed nuisance, the outcome distribution conditional on its response history is iid Bernoulli$(p)$. For $n>0$ the center and endpoint risks of $q_\alpha$ are
\begin{equation}\label{eq:ae}
 A_\alpha(n)=\frac{n}{4(n+2\alpha)^2},\qquad
 E_\alpha(n)=\frac{\alpha^2}{(n+2\alpha)^2},
\end{equation}
with $A_\alpha(0)=0$, $E_\alpha(0)=1/4$. The whole-parameter risk is
\begin{equation}\label{eq:quadraticrisk}
 A+4(E-A)(p-1/2)^2,
\end{equation}
after averaging these coefficients over $N$. Its symmetric Beta average is \eqref{eq:bcount}, the minimum terminal Bayes risk. Proposition~\ref{prop:bayes} implies that $F_R(\alpha)$ is also the minimum Bayes risk over all full procedures under $\Beta(\alpha,\alpha)\times\nu$. Hence $V_F\ge\sup_{\alpha>0}F_R(\alpha)$. The corresponding static lower bound holds as well. Taking $\alpha\downarrow0$ gives the proper half-half endpoint prior lower bound.

For a fixed tree put $f_\pi(\alpha)=\E_\nu^\pi b_\alpha(N)$. Symmetrization preserves $f_\pi$ and equalizes its risk across nuisance states. Denote the averaged coefficients by $A_\pi,E_\pi$. Direct differentiation gives
\begin{equation}\label{eq:slope}
 f_\pi'(\alpha)=-\frac{2(E_\pi-A_\pi)}{(2\alpha+1)^2}.
\end{equation}
For $n>0$, the derivative of $b_\alpha(n)$ is
$(n-4\alpha^2)/[2(2\alpha+1)^2(n+2\alpha)^2]$; for $n=0$ it is $-1/[2(2\alpha+1)^2]$. These identities also hold as right limits at zero, where $A_0=\E[1/(4N);N>0]$ and $E_0=\Pp(N=0)/4$.

The minimum of finitely many such continuous functions is continuous. It is positive at positive $\alpha$ and bounded by $1/[4(2\alpha+1)]$, hence tends to zero at infinity and attains its maximum at finite $\alpha$. At an interior maximum the left derivative of this finite minimum is the largest slope of its active functions and the right derivative is the smallest. They straddle zero. A mixture of at most two active trees makes the average slope zero. By \eqref{eq:slope}, its averaged $A,E$ coincide. Symmetrization and \eqref{eq:quadraticrisk} now give a risk constant in every parameter, equal to its Bayes value $F_R(\alpha_*)$. This matches the arbitrary-full lower bound.

At a maximum at zero, the right derivative is nonpositive, so some active tree has $E_0\ge A_0$. Its worst risk is $E_0=F_R(0)$, matching the endpoint lower bound. Applying the same argument to static schedules proves the second formula. All randomizations are independent seeds known to the estimator and are averaged before taking worst risk.
\end{proof}

\begin{remark}
A symmetrized full policy need not have ancillary $N$. With two equal-cost contacts, choose the first channel uniformly and repeat it only after a response with $Y=1$, stopping otherwise. Then $\Pp(N=2)=p(\rho_1^2+\rho_2^2)/2$. Group averaging cannot replace Proposition~\ref{prop:bayes}. Likewise, that proposition cannot replace the slope-balancing saddle construction. With several nuisance orbits, symmetry only equalizes risks within each orbit; a further joint saddle is required.
\end{remark}

\section{A sharp two-contact classification}\label{sec:two}
Consider $c=(1,1)$, $B=2$ and $\Gamma=\{(a,b),(b,a)\}$, $a,b\in[0,1]$. Write $\Delta=(a-b)^2/2$. A balanced static schedule contacting each channel once gives
\[
 f_S=((1-a)(1-b),\ a+b-2ab,\ ab).
\]
A uniformly labeled repeated-channel static schedule gives
\[
 f_D=(f_{S,0}+\Delta,\ f_{S,1}-2\Delta,\ f_{S,2}+\Delta).
\]
Choose the first channel uniformly, repeat it after a response and switch after a failure. In either nuisance state its count law is
\[
 f_R=(f_{S,0},\ f_{S,1}-\Delta,\ f_{S,2}+\Delta).
\]

\begin{corollary}\label{cor:two}
In this experiment,
\begin{equation}\label{eq:twovalue}
 V_S=\max_{\alpha\ge0}\sum_{n=0}^2f_{S,n}b_\alpha(n),\quad
 V_F=V_R=\max_{\alpha\ge0}\sum_{n=0}^2f_{R,n}b_\alpha(n).
\end{equation}
Moreover,
\begin{equation}\label{eq:iff}
 V_R<V_S\quad\Longleftrightarrow\quad
 a\ne b\ \hbox{ and }\ 4(a+b)-5ab>2.
\end{equation}
These comparisons include all pre-data static randomization, all full policies and all terminal estimators.
\end{corollary}
\begin{proof}
For a symmetric proper Beta prior let $B_\alpha(f)=\sum f_nb_\alpha(n)$. Exact subtraction gives
\begin{align}\label{eq:bayesdiff}
 B_\alpha(f_D)-B_\alpha(f_S)
   &=\frac{(a-b)^2}{8(\alpha+1)(2\alpha+1)^2},\\
 B_\alpha(f_S)-B_\alpha(f_R)
   &=\frac{\alpha(a-b)^2}{8(\alpha+1)(2\alpha+1)^2}.\nonumber
\end{align}
Thus balanced sampling is Bayes optimal among all static schedules for every $\alpha>0$. A schedule using fewer contacts cannot improve on an extension whose estimator can ignore the added data.

Under the uniform nuisance prior, after a first response the predictive response probability of repeating exceeds that of switching by $(a-b)^2/(a+b)$; after a failure the predictive advantage of switching is $(a-b)^2/(2-a-b)$. A zero denominator indicates an impossible branch. By Proposition~\ref{prop:bayes}, the choice need not use its $Y$ value. Hence the rule giving $f_R$ is Bayes optimal even among full policies for every proper symmetric Beta prior. Theorem~\ref{thm:orbit} proves \eqref{eq:twovalue}.

For completeness the global optimizing alternative can be checked without an inference from a derivative at zero alone. For a fixed count law let $A_\alpha(f),E_\alpha(f)$ average \eqref{eq:ae}. Unless $f_0=1$, $A$ strictly decreases, $E$ increases, $A_\infty=0$, $E_\infty=1/4$. If $A_0\le E_0$, the sample mean with half at zero has worst risk $E_0$, matched by the endpoint prior. Otherwise there is a unique positive crossing $A_\alpha=E_\alpha$, and the corresponding shrinkage rule has constant risk matched by its Beta prior. The static endpoint condition is
\[
 f_{S,0}\ge f_{S,1}+f_{S,2}/2
 \quad\Longleftrightarrow\quad4(a+b)-5ab\le2.
\]
In that region $A_0(f_R)\le A_0(f_S)\le E_0$ and both values are $f_{S,0}/4$. Outside it $V_S>f_{S,0}/4$. If the $f_R$ optimizer is zero this gives the strict inequality. If it is positive and $a\ne b$, the second identity in \eqref{eq:bayesdiff} gives $V_S\ge B_{\alpha_R}(f_S)>B_{\alpha_R}(f_R)=V_R$. When $a=b$ the laws coincide. This includes rates 0 or 1 and the equality boundary; for $a=b=0$ the common value is $1/4$.
\end{proof}

The count improvement transfers probability $\Delta$ from one to two responses without changing the all-failure probability. When the endpoint prior controls the worst risk, this improvement has no minimax value, although every proper symmetric Beta prior benefits when $a\ne b$.

\section{An exact cost-asymmetric saddle}\label{sec:asym}
We now take $c=(1,2)$, $B=5$, $r=1/2000$ and two nuisance states $\rho^0=(1-r,r)$, $\rho^1=(r,1-r)$. This experiment violates the exchange-cost hypothesis of Theorem~\ref{thm:orbit}. Its solution instead uses a parameter-correlated finite prior and a certified isolated root. With unequal costs the two channels are no longer interchangeable, so group averaging cannot equalize the two nuisance states. One might then expect outcome feedback to help, for example by deciding whether to buy an expensive contact only after seeing an outcome. The example shows that it does not: an exact saddle is again attained by a response-only rule. The price is that the least-favourable prior must couple $p$ with the nuisance state (Proposition~\ref{prop:nonproduct}), and the saddle can only be located numerically, so we certify it by exact interval arithmetic.

Call the cheap and expensive channels $C,E$. Two response trees first contact $C$ and $E$ once. At response pair $(1,0)$, tree $T_C$ contacts $C$ twice more, while $T_E$ contacts $E$ once more. At any other response pair both trees contact $E$ once more. Each path costs five and has at most three responses. Choose $T_C$ with pre-data probability $\theta$ and $T_E$ otherwise. The estimator below is Bayes for the same prior under either tree.

Let $x=(L,u_0,v_0,u_1,v_1,\theta)$, $\lambda=L/10^7$. Give nuisance 0 mass $\lambda$ and nuisance 1 mass $1-\lambda$. Conditional on nuisance $j$, use a symmetric prior on $\{0,1/5,1/2,4/5,1\}$ with moments $u_j=\E_j[p(1-p)]$, $v_j=\E_j[p^2(1-p)^2]$. Its edge, middle and central weights are
\begin{equation}\label{eq:weights}
 w_{j,H}=\frac{v_j-(4/25)u_j}{9/400},\qquad
 w_{j,M}=\frac{u_j-w_{j,H}/4}{4/25},\qquad
 w_{j,E}=1-w_{j,H}-w_{j,M}.
\end{equation}
The five weights are $w_{j,E}/2,w_{j,M}/2,w_{j,H},w_{j,M}/2,w_{j,E}/2$.

For response statistics $z=(a_C,b_C,a_E,b_E)$ define
\begin{align}\label{eq:momentmix}
 \ell_j(z)&=(\rho_C^j)^{a_C}(1-\rho_C^j)^{b_C}
                (\rho_E^j)^{a_E}(1-\rho_E^j)^{b_E},\nonumber\\
 m_z&=\lambda\ell_0(z)+(1-\lambda)\ell_1(z),\\
 U_z&=\frac{\lambda\ell_0(z)u_0+(1-\lambda)\ell_1(z)u_1}{m_z},\quad
 Q_z=\frac{\lambda\ell_0(z)v_0+(1-\lambda)\ell_1(z)v_1}{m_z}.\nonumber
\end{align}
Write $n=a_C+a_E$ and let $s$ count observed successes. Direct atomic integration gives the posterior mean $d_{n,s}(U,Q)$:
\begin{center}\small
\begin{tabular}{cl}\toprule
$n$&$(d_{n,0},\ldots,d_{n,n})$\\\midrule
0&$(1/2)$\\
1&$(2U,1-2U)$\\
2&$\bigl(U/(1-2U),1/2,1-U/(1-2U)\bigr)$\\
3&$\bigl((U-2Q)/(1-3U),2Q/U,1-2Q/U,1-(U-2Q)/(1-3U)\bigr)$\\\bottomrule
\end{tabular}
\end{center}
Use $d_{n,s}(U_z,Q_z)$ as the actual terminal estimate.

Here is a finite exact specification of $x$. For $T_E$, include one path $z=(c,1-c,e+j,2-e-j)$ for each $c,e,j\in\{0,1\}$. For $T_C$ retain those paths except when $(c,e)=(1,0)$; replace that pair by $z=(1+k,2-k,0,1)$, $k=0,1,2$, with multiplicities $\binom2k$. Identical response statistics from different paths are added. Let these multisets be $\mathcal L_C,\mathcal L_E$. Define
\begin{equation}\label{eq:poly}
\begin{split}
 \mathcal R_j(p;x)={}&\sum_{t=C,E}w_t\sum_{(z,a)\in\mathcal L_t}a\ell_j(z)\\
 &\qquad\times\sum_{s=0}^n\binom ns p^s(1-p)^{n-s}
                [d_{n,s}(U_z,Q_z)-p]^2,
\end{split}
\end{equation}
where $w_C=\theta,w_E=1-\theta$. Its coefficients are rational functions of $x$ and its degree is at most five.

Put $h_{n,s}(U,Q)=\E[p^s(1-p)^{n-s}]$ for a symmetric prior with these moments. The vectors for $n=0,1,2,3$ are
\begin{align*}
 &(1),\quad(1/2,1/2),\quad(1/2-U,U,1/2-U),\\
 &(1/2-3U/2,U/2,U/2,1/2-3U/2).
\end{align*}
Its averaged posterior risk is
\[
 g_n(U,Q)=1/2-U-\sum_{s=0}^n\binom ns h_{n,s}(U,Q)d_{n,s}(U,Q)^2.
\]
Set $z_k^C=(1+k,2-k,0,1)$ and $z_k^E=(1,0,k,2-k)$, and define
\begin{equation}\label{eq:tie}
\begin{split}
 D(x)={}&\sum_{k=0}^2\binom2k m_{z_k^C}g_{1+k}(U_{z_k^C},Q_{z_k^C})\\
 &-\sum_{k=0}^1m_{z_k^E}g_{1+k}(U_{z_k^E},Q_{z_k^E}).
\end{split}
\end{equation}
The six defining equations are
\begin{equation}\label{eq:system}
 F(x)=\begin{pmatrix}
 [p^2]\mathcal R_0\\ [p^4]\mathcal R_0\\
 [p^2]\mathcal R_1\\ [p^4]\mathcal R_1\\
 \mathcal R_0(0;x)-\mathcal R_1(0;x)\\ 10^7D(x)
 \end{pmatrix}=0.
\end{equation}
The first five equations equalize risk; the sixth makes the two actions indifferent at the randomization branch. An explicit rational isolation box is given in Appendix~\ref{app:box}.

\begin{theorem}[Unequal-cost saddle]\label{thm:asym}
System~\eqref{eq:system} has a unique zero $x_*$ in the box of Appendix~\ref{app:box}. Its prior weights are strictly positive and $0<\lambda_*,\theta_*<1$. In the experiment of this section,
\begin{equation}\label{eq:asymvalue}
 V_F=V_R=V_*:=\mathcal R_0(0;x_*).
\end{equation}
The response-only mixture and terminal estimate above attain this value at every $p$ in both nuisance states. The ten-point prior is least favourable for the entire full-feedback class. In particular,
\[
\begin{gathered}
 V_*>0.0429061342528548909317426699137608697654,\\
 V_*<0.0429061342528548909317426699137608697659,
\end{gathered}
\]
and $\theta_*=0.6497480330156957348346\ldots$. The exact value is defined by \eqref{eq:system} together with its box.
\end{theorem}
\begin{proof}
We first describe how the exact finite certificate proves existence and Bayes optimality, then the analytic upper and lower bounds. The rational certificate supplies a matrix $M$ and outward Jacobian bounds on the box $\mathcal B$ with
\begin{equation}\label{eq:contraction}
\begin{gathered}
 \norm{I-MJ_F(\mathcal B)}_\infty<1.3\times10^{-34},\\
 x_c-MF(x_c)+(I-MJ_F(\mathcal B))(\mathcal B-x_c)
       \subset\operatorname{int}\mathcal B.
\end{gathered}
\end{equation}
All denominators in $F$ are positive there. The mean-value formula implies that $x\mapsto x-MF(x)$ maps the closed box into itself and is a contraction. Since $MJ_F(x_c)$ is invertible by the first inequality, $M$ is invertible. Thus its unique fixed point is a unique zero of $F$ in the box. Positivity of \eqref{eq:weights}, $\lambda$, $1-\lambda$ and both action probabilities is verified on the same box.

For this prior, define the unnormalized moments and terminal loss
\begin{equation}\label{eq:unnorm}
 M_k(z,s)=\sum_{(j,p,w)}w\ell_j(z)p^{s+k}(1-p)^{n-s},\qquad
 H(z,s)=M_2(z,s)-\frac{M_1(z,s)^2}{M_0(z,s)}.
\end{equation}
All $M_0$ are positive. The minimum unnormalized continuation risk over full policies is the finite recursion
\begin{align}\label{eq:fullbellman}
 W(e,z,s)&=\min\{H(z,s),\min_{c_h\le e} J_h(e,z,s)\},\\
 J_h(e,z,s)&=W(e-c_h,z+h0,s)+W(e-c_h,z+h1,s)\nonumber\\
 &\hspace{27mm}+W(e-c_h,z+h1,s+1).\nonumber
\end{align}
The three children correspond to failure, responding outcome 0 and responding outcome 1. Likelihood masses are already included, so no additional predictive-probability multiplier belongs in this unnormalized recursion. Posterior moments are sufficient at every full history by \eqref{eq:factor}; parameter-free choices do not resample nuisance. Posterior means minimize terminal risk, and randomization cannot improve the minimum. Hence \eqref{eq:fullbellman} includes arbitrary full allocation, stopping and estimation.

The supplied finite certificate checks every affordable action at every sufficient-statistic state reachable under any action. It represents each continuation as an integer linear combination of the losses in \eqref{eq:unnorm}. Symmetry gives $H(z,s)=H(z,n-s)$. At the critical one-response prefix the full-history action difference for either first outcome is $D/2$, identically, so \eqref{eq:system} certifies the tie. Other exact ties arise from the same terminal path multiset or a rational multiple of this difference. All remaining action differences are strictly positive by outward rational intervals, with minimum recorded lower bound greater than $4.62\times10^{-14}$. These identities and signs show that the optimal root continuation is equivalent to $T_C$, and that $T_E$ has the same Bayes value. Their mixture is full-Bayes optimal. Lemma~\ref{lem:law} makes this a lower bound for arbitrary measurable randomized procedures. Appendix~\ref{app:certificate} specifies the independently implemented proof consumer and the complete finite data, rather than a floating-point optimum.

For the attaining response-only design, \eqref{eq:poly} is the exact frequentist risk. To see this directly, fix either tree. Its response paths have the multisets listed above: only the prefix $(c,e)=(1,0)$ changes, and two further cheap responses have multiplicities $1,2,1$. Conditional on any such response path, allocation ignores outcomes and its $n$ responding observations are independent Bernoulli variables with parameter $p$. Grouping their $2^n$ outcome strings by their success count gives the binomial coefficient in \eqref{eq:poly}. The initial independent tree choice then gives its two mixture weights. Thus the compressed formula is the risk of a feasible rule, not merely a formal polynomial. Its estimate satisfies $d_{n,n-s}=1-d_{n,s}$, hence $\mathcal R_j(p;x)=\mathcal R_j(1-p;x)$. A symmetric polynomial of degree at most five is a linear combination of $1,p(1-p),p^2(1-p)^2$. Its $p^4$ and $p^2$ coefficients vanishing makes it constant. The first five equations in \eqref{eq:system} thus give
\[
 \mathcal R_0(p;x_*)=\mathcal R_1(p;x_*)=V_*,\quad0\le p\le1.
\]
Its prior-average risk is $V_*$ and is the full Bayes optimum just established. Every full rule has worst risk at least that Bayes optimum. The response-only rule has worst risk $V_*$, proving \eqref{eq:asymvalue}. The displayed enclosure follows by interval evaluation on the isolating box, not by comparison of two approximately equal optimization values.
\end{proof}

\subsection{Why a product prior cannot give this saddle}
\begin{proposition}\label{prop:nonproduct}
No independent prior $G(dp)\times\nu(dj)$ is least favourable for the experiment in Section~\ref{sec:asym}. Indeed, its optimized full Bayes lower bounds satisfy
\[
 \sup_{G,\nu}\inf_{D\in\Pi_{\mathrm{full}}}\int r(D;p,j)\,d(G\times\nu)
                 \le V_*-4\times10^{-17}.
\]
Hence a search restricted to product priors cannot find the saddle, even though the minimax values with and without outcome feedback coincide.
\end{proposition}
\begin{proof}
The attaining rule's terminal response states $z_A=(0,1,1,1)$ and $z_B=(1,0,0,2)$, both with $n=1,s=0$, have estimates $d_A=2U_{z_A}$ and $d_B=2U_{z_B}$. The box verifies $u_1>u_0$, $0<\lambda<1$, and their nuisance likelihood ratios are $r/(1-r)$ and $((1-r)/r)^3$. Thus $d_A>d_B$. Both histories have positive response probability in either state; $z_B$ occurs on a branch of probability $1-\theta_*>0$.

Under any independent prior other than $G=\delta_1$, their unique posterior mean is the same number $\int p(1-p)\,dG/\int(1-p)\,dG$. If that prior were least favourable, the equalizer, with average risk $V_*$, would be Bayes optimal and would have to use this mean at both histories, a contradiction. The excluded prior has Bayes risk zero.

For the quantitative bound put $e=\E_G(1-p)$ and $\delta=d_A-d_B$. Retain one response path for $z_A$ and the expensive continuation path for $z_B$, each with its observed zero. Their prior probabilities are at least $e r^2(1-r)$ and $e(1-\theta_*)r^3$. Replacing both estimates by their common posterior mean improves average risk by at least
\[
 \frac{w_Aw_B}{w_A+w_B}\delta^2
        \ge\frac{e(1-\theta_*)r^3\delta^2}{2}.
\]
If $e\ge V_*/2$, use the verified inequalities $V_*>0.04$, $1-\theta_*>0.35$, $\delta>0.01$; the improvement exceeds $4.375\times10^{-17}$. If $e<V_*/2$, the constant estimate 1 has Bayes risk $\E_G(1-p)^2\le e<V_*/2$, already smaller than the asserted bound. This proves the statement for every $G,\nu$, including degeneracies.
\end{proof}

\subsection{Admissibility and nonuniqueness}
\begin{proposition}\label{prop:admiss}
The attaining response-only rule is admissible within the full class. Every minimax rule has risk $V_*$ at every parameter. Minimax allocation is not unique.
\end{proposition}
\begin{proof}
A legal full rule's risk is a polynomial of degree at most seven in $p$, for fixed nuisance: each complete path uses at most five observed outcomes, its probability has degree at most five, and squared loss adds at most two. Averaging arbitrary independent randomization preserves this degree bound. For a proposed dominator, $V_*-r(D;p,j)$ is nonnegative on $[0,1]$. Least favourability and optimality force its integral under the positive ten-point prior to be zero. It vanishes at all five support points in each nuisance state. The three interior zeros have even multiplicity at least two; the two endpoint zeros each have multiplicity at least one. A nonzero polynomial would require degree at least eight. It is therefore identically zero, excluding strict domination. The argument also applies to every minimax rule. Finally, either order of the first cheap/expensive contacts gives the same response-statistic distribution, subsequent decision and risk, while preserving budget; these are distinct attaining allocation trees.
\end{proof}

Response feedback does have value in this asymmetric example. A common static lower-bound prior takes $p\sim\Beta(1/\sqrt2,1/\sqrt2)$ independently of nuisance masses $(1/3,2/3)$. For any preassigned counts $(n_C,n_E)$, $n_C+2n_E\le5$, its response-count probability generating function is
\begin{align*}
 G_{n_C,n_E}(t)={}&\tfrac13[r+(1-r)t]^{n_C}[1-r+rt]^{n_E}\\
 &+\tfrac23[1-r+rt]^{n_C}[r+(1-r)t]^{n_E}.
\end{align*}
The arbitrary-estimator Bayes risk is $\sum_m[t^m]G_{n_C,n_E}(t)b_{1/\sqrt2}(m)$. Exhausting the twelve legal count pairs gives its minimum at $(1,2)$,
\[
 L_S=\frac{9299030989}{56000000000}
          -\frac{18568087969}{224000000000}\sqrt2
       >0.0488253666.
\]
Every static randomization has Bayes risk at least this common-prior minimum, so $V_S\ge L_S$ and $V_S-V_R>0.0059192323$. We do not compute $V_S$ exactly here. The supplement checks the radical inequality by rational squaring and enumerates all schedules.

\section{Bayes feedback need not be valueless}\label{sec:counter}
\begin{proposition}\label{prop:counter}
In the two-contact equal-cost experiment with rates $(3/4,1/4)$ and their swap, give prior mass $1/2$ to $(p,\rho)=(1/5,(3/4,1/4))$ and $(4/5,(1/4,3/4))$. Then
\[
 b_F=\frac{23082273}{629000000},\qquad
 b_R=\frac{12195549}{327080000},
\]
\[
 b_R-b_F=\frac{16281}{27625000}>0.
\]
Its frequentist full and response-only minimax values are nevertheless equal.
\end{proposition}
\begin{proof}
Condition on first channel 1. Enumerating the three first observations gives the following unnormalized posterior squared-risk contributions for the second-channel choice. In each cell terminal estimates have already been optimized by posterior means.
\begin{center}
\begin{tabular}{ccc}\toprule
First observation&Next channel 1&Next channel 2\\\midrule
$R=0$&$837/40000$&$19467/769600$\\
$R=1,Y=0$&$10827/1450000$&$7749/1258000$\\
$R=1,Y=1$&$2403/250000$&$4509/442000$\\\bottomrule
\end{tabular}
\end{center}
The full optimum sums row minima. Response-only allocation must use one common next action for the last two rows, so it minimizes their column sums instead. The first channel 2 gives the same values by simultaneous label/outcome reflection. Independent randomized actions cannot improve these minima, and additional data can be ignored, so stopping cannot do better. Each entry follows by summing $M_2-M_1^2/M_0$ over its terminal outcomes. The resulting fractions are those displayed. Theorem~\ref{thm:orbit} gives the minimax equality for this single cost-preserving orbit.
\end{proof}

\section{Discussion}\label{sec:discussion}
The results distinguish three statements often conflated in adaptive design: an optimal design for a fixed prior, a worst-risk equality between information classes, and comparison of experiments across decision problems. Product conjugacy supplies the first. A feasible equalizer matched to an arbitrary-full-policy Bayes lower bound supplies the second. Neither supplies the third. The correlated-prior example is a direct obstruction to an all-prior interpretation.

The statistical restrictions are material. All channels observe the same Bernoulli outcome distribution; response is noninformative conditional on fixed parameters; costs are known and result-independent. Informative nonresponse, channel-specific means, changing response states, expected rather than pathwise budgets and sequential target decisions require different arguments. The exact asymmetric example proves that cost asymmetry does not force outcome feedback value, not that it never creates it. Its nonproduct least-favourable prior shows that searching only independent priors can fail even when the final minimax allocation ignores outcomes.

The main contribution is the exact minimax saddle for controlled experiments and its consequences for the value of feedback; random-count risks, group averaging and the conjugate factorization are standard ingredients. The asymmetric instance supports and delimits the general result. Extending the analysis to several nuisance orbits, informative nonresponse and asymptotic budgets are natural directions for future work.

\appendix

\section{Exact root box and certificate interpretation}\label{app:box}
The box has componentwise radius $10^{-40}$ about the following exact decimal rational centers, ordered as $(L,u_0,v_0,u_1,v_1,\theta)$:
\begin{center}\footnotesize\ttfamily
\begin{tabular}{l}
1.5546796686568088491498779422143555056830072715449037624796973752\\
0.14126992023989189784295197340241505233685476496627090421053190865\\
0.027102721330304065718853360459422214213584016742722069357652737615\\
0.14642079608220506620748713166079454511069532453849214547088378506\\
0.028311493143517131572331703537711033968392283182344633552526068618\\
0.64974803301569573483460409103551566648191858333485940970284546116
\end{tabular}
\end{center}
The complete rational matrix, outward Jacobian enclosures and contraction image are stored in the machine-readable supplementary certificate. The root specification \eqref{eq:system} is independent of a floating-point solver: supplied rational data satisfy \eqref{eq:contraction} using integer arithmetic. Clearing the nonvanishing denominators yields a nonsingular isolated polynomial system over $\mathbb Q$, an algebraic definition of the value; no expression by radicals is required.

\section{Certified computation and reproducibility}\label{app:certificate}
The supplementary proof consumer uses only Python's standard library. It reconstructs the prior from its ten atoms, posterior means from moment ratios, and frequentist risk polynomials from individual observation paths. It independently implements outward binary fixed-point interval arithmetic and differentiated interval operations. It reads the rational matrix and root box as proposed proof data and verifies the contraction and self-map, without importing the original generator or its arithmetic core. It then enumerates all 127 sufficient-statistic states of \eqref{eq:fullbellman}; 89 nonidentical comparisons have strictly positive interval margins and nine comparisons are exact ties. All are finite proofs, not parameter-grid maxima.

For ties the consumer expresses a continuation as an integer combination of terminal losses. Identifying $s$ with $n-s$ uses the exact reflection symmetry of the atomic prior. A difference proportional to the critical $D/2$ expression is zero at the isolated root by \eqref{eq:system}. An interval merely containing zero is never accepted as proof of equality. All other differences are certified by strict signs. The full action space includes stop and every affordable channel; no off-policy branch is dropped.

The ancillary files \texttt{anc/verify.py} and \texttt{anc/certificate.json} reproduce these checks using Python~3 and its standard library. Execute
\begin{quote}\small\ttfamily
cd anc \\
python verify.py --certificate certificate.json
\end{quote}
Exact fractions verify Proposition~\ref{prop:counter}, and radical interval bounds verify the static lower bound. The consumer is independent of the program that generated the certificate. The probability model and minimax theorems have not been formalized in a proof assistant; the finite arithmetic described here is the only computational component of the proofs.

\section*{Ancillary files}
The ancillary files of this arXiv submission contain the standard-library Python consumer \texttt{anc/verify.py}, the exact rational certificate \texttt{anc/certificate.json} and reproduction instructions. They check the local root isolation, all full-policy Bellman states, strict comparisons and exact ties, the displayed risk enclosure, the twelve static schedules and the correlated-prior example of Section~\ref{sec:counter}. The finite arithmetic checks support the analytic proofs rather than replacing their statistical interpretation.

\section*{Statement on the use of AI tools}
During the preparation of this work the author used OpenAI Codex (October 2026) to assist with mathematical derivations, literature retrieval, proof checking, development of the verification code, and drafting of the manuscript text, and Anthropic's Claude (Claude Opus 5.5, November 2026) for an independent computational re-check of the main results and for formatting the manuscript. After using these tools, the author reviewed and edited the content as needed and takes full responsibility for the content of the published article.

\section*{Declaration of competing interest}
The author declares no competing interests.

\section*{Funding}
This research did not receive any specific grant from funding agencies in the public, commercial, or not-for-profit sectors.

\bibliographystyle{plainnat}
\bibliography{references}
\end{document}